\documentclass{amsart}

\usepackage[english]{babel}

\newtheorem{theorem}{Theorem}[section] 

\newcommand{\C}{{\mathbb C}}       \newcommand{\B}{{\mathbb B}} 
    
 \DeclareMathOperator{\Aut}{Aut}

\begin{document}


\title{Some Irreducible Representations of the Braid Group $\B_n$ of Dimension greater than $n$}

\author{Claudia Mar\'{i}a Egea, Esther Galina}

\address{Facultad de Matem\'atica Astronom\'{i}a y F\'{i}sica, Universidad Nacional de C\'ordoba, C\'ordoba, Argentina} \email{
cegea@mate.uncor.edu, galina@mate.uncor.edu}

\begin{abstract} For any $n\geq 3$, we construct a family of finite dimensional irreducible representations of the braid group $\B_n$. Moreover,
we give necessary conditions for a member of this family to be irreducible. In particular we give a explicitly irreducible subfamily $(\phi_m,
V_m)$, $1\leq m<n$, where $\dim V_m=\left( \begin{smallmatrix} n\\\noalign{\medskip}m \end{smallmatrix}\right)$. The representation obtained in
the case $m=1$ is equivalent to the standard representation. \end{abstract}


 \thanks{This work was partially supported by CONICET, SECYT-UNC, FONCYT.} \maketitle

\section{Introduction}\label{Introduccion}

The braid group of $n$ strings $\B_n$, is defined by generators and relations as follows $$ \B_{n}=<\tau_{1}, \dots,\tau_{n-1}>_{/\sim} $$

$$ \sim=\{ \tau_k\tau_j=\tau_j\tau_k, \textrm{ if }|k- j|>1;  \ \
        \tau_k\tau_{k+1}\tau_k=\tau_{k+1}\tau_{k}\tau_{k+1} \ \ 1\leq k\leq n-2\ \}
$$

We will consider finite dimensional complex representations of $\B_n$; that is pairs $(\phi,V)$ where $$ \phi : \B_n\rightarrow \Aut(V) $$ is a
morphism of groups and $V$ is a complex vector space of finite dimension.

In this paper, we will construct a family of finite dimensional complex representations of $\B_n$ that contains the standard representations.
Moreover, we will give necessary conditions for a member of this family to be irreducible. In this way, we can find explicit families of
irreducible representations. In particular, we will define a subfamily of irreducible representations $(\phi_m, V_m)$, $1\leq m<n$, where $\dim
V_m=\left( \begin{smallmatrix} n\\\noalign{\medskip}m \end{smallmatrix}\right)$ and the corank of $\phi_m$ is equal to
$\frac{2(n-2)!}{(m-1)!(n-m-1)!}$.

This family of representations can be useful in the progress of classification of the irreducible representations of $\B_n$. As long as we
known, there are only few contributions in this sense, some known results are the following ones. Formanek classified all the irreducible
representations of $\B_n$ of dimension lower than $n$ \cite{F}. Sysoeva did it for dimension equal to $n$ \cite{S}. Larsen and Rowell gave some
results for unitary representation of $\B_n$ of dimension multiples of $n$. In particular, they prove there are not irreducible representations
of dimension $n+1$. Levaillant proved when the Lawrence-Krammer representation is irreducible and when it is reducible \cite{L}.

\section{Construction and Principal Theorems}

In this section, we will construct a family of representations of $\B_n$ that we believe to be new, and we will obtain a subfamily of
irreducible representations.

We choose $n$ non negative integers $z_1, z_2, \dots, z_n$, not necessarily different. Let $X$ be the set of all the possible $n$-tuples
obtained by permutation of the coordinates of the fixed $n$-tuple $(z_1, z_2, \dots, z_n)$. For example, if the $z_i$ are all different, then
the cardinality of $X$ is $n!$. Explicitly, if $n=3$, $$ X=\{(z_1, z_2, z_3), (z_1, z_3, z_2), (z_2, z_1, z_3), (z_2, z_3, z_1), (z_3, z_1,
z_2), (z_3, z_2, z_1)\} $$

Or if $z_1=z_2=1$ and $z_i=0$ for all $i=3, \dots, n$, then the cardinality of $X$ is $\left( \begin{smallmatrix} n\\\noalign{\medskip}2
\end{smallmatrix}\right)=\frac{n(n-1)}{2}$. Explicitly, for $n=3$ $$ X=\{(1,1,0), (1,0,1), (0,1,1)\} $$

Let $V$ be a complex vector space with orthonormal basis $\beta=\{v_x : x\in X\}$. Then the dimension of $V$ is the cardinality of $X$.

We define $\phi:\B_n \rightarrow \Aut (V)$, such that $$ \phi(\tau_k)(v_x)= q_{x_k, x_{k+1}} v_{\sigma_k(x)} $$ where $q_{x_k, x_{k+1}}$ is a
non-zero complex number that depends on $x=(x_1, \dots, x_n)$, but, it only depends on the places $k$ and $k+1$ of $x$; and $$ \sigma_k(x_1,
\dots, x_n)=(x_1, \dots,x_{k-1}, x_{k+1}, x_k, x_{k+2}, \dots, x_n) $$ With this notations, we have the following theorem,

\begin{theorem} $(\phi, V)$ is a representation of the braid group $\B_n$. \end{theorem} \begin{proof} We need to check that $\phi(\tau_k)$
satisfy the relations of the braid group. We have for $j\neq k-1, k, k+1$ that $$ \phi(\tau_k)\phi(\tau_j)(v_x)= \phi(\tau_k)(q_{x_j, x_{j+1}}
v_{\sigma_j(x)})=q_{x_j, x_{j+1}} q_{x_k, x_{k+1}} v_{\sigma_k \sigma_j(x)} $$ On the other hand $$ \phi(\tau_j)\phi(\tau_k)(v_x)=
\phi(\tau_j)(q_{x_k, x_{k+1}} v_{\sigma_k(x)})=q_{x_k, x_{k+1}} q_{x_j, x_{j+1}} v_{\sigma_j \sigma_k(x)} $$ As $\sigma_k \sigma_j(x)=\sigma_j
\sigma_k(x)$, if $|j-k|>1$, then $\phi(\tau_k)\phi(\tau_j)=\phi(\tau_k)\phi(\tau_j)$ if $|j-k|>1$.

In the same way, we have $$ \begin{aligned} \phi(\tau_k)\phi(\tau_{k+1})\phi(\tau_k)(v_x)&= \phi(\tau_k)\phi(\tau_{k+1})(q_{x_k, x_{k+1}}
v_{\sigma_k(x)})\\ &=\phi(\tau_k)(q_{x_k, x_{k+1}} q_{x_{k}, x_{k+2}} v_{\sigma_{k+1} \sigma_k(x)})\\ &= q_{x_k, x_{k+1}} q_{x_{k}, x_{k+2}}
q_{x_{k+1}, x_{k+2}} v_{\sigma_k \sigma_{k+1} \sigma_k(x)} \end{aligned} $$ Similarly, $$ \begin{aligned}
\phi(\tau_{k+1})\phi(\tau_k)\phi(\tau_{k+1})(v_x)&= \phi(\tau_{k+1})\phi(\tau_k)(q_{x_{k+1}, x_{k+2}} v_{\sigma_{k+1}(x)})\\
&=\phi(\tau_{k+1})(q_{x_{k+1}, x_{k+2}} q_{x_{k}, x_{k+2}} v_{\sigma_k \sigma_{k+1}(x)})\\ &= q_{x_{k+1}, x_{k+2}} q_{x_{k}, x_{k+2}} q_{x_k,
x_{k+1}} v_{\sigma_{k+1} \sigma_k \sigma_{k+1}(x)} \end{aligned} $$ As $\sigma_k \sigma_{k+1} \sigma_k(x)=\sigma_{k+1} \sigma_k
\sigma_{k+1}(x)$, for all $k$ and $x\in X$, then $\phi(\tau_k)\phi(\tau_{k+1})\phi(\tau_k)=\phi(\tau_{k+1})\phi(\tau_k)\phi(\tau_{k+1})$ for all
$k$. \end{proof}

As $\beta$ is an orthonormal basis, we have that, $$ <\phi(\tau_k)v_y, v_x> = <q_{y_k, y_{k+1}}v_{\sigma_k(y)}, v_x> = <v_y,
\overline{q_{x_{k+1}, x_k}} v_{\sigma_k(x)}> $$ then, $$ (\phi(\tau_k))^*(v_x)= \overline{q_{x_{k+1}, x_k}} v_{\sigma_k(x)} $$ therefore,
$\phi(\tau_k)$ is self-adjoint if and only if $q_{x_{k+1}, x_k}=\overline{q_{x_k, x_{k+1}}}$ for all $x\in X$. In particular, if $x_k=x_{k+1}$
then $q_{x_k, x_{k+1}}$ is a real number. In the same way, $\phi(\tau_k)$ is unitary if and only if $|q_{x_k, x_{k+1}}|^2=1$ for all $x\in X$.

Now, we will give a subfamily of irreducible representations.

\begin{theorem} \label{irreducible} If $\phi(\tau_k)$ is a self-adjoint operator for all $k$, and for any pair $x,y\in X$, there exists $j$,
$1\leq j\leq n-1$, such that $|q_{x_j, x_{j+1}}|^2\neq |q_{y_j, y_{j+1}}|^2$, then $(\phi, V)$ is an irreducible representation of the braid
group $\B_n$. \end{theorem} \begin{proof} Let $W\subset V$ be a non-zero invariant subspace. It is enough to prove that $W$ contains one of the
basis vectors $v_x$. Indeed, given $y\in X$, there exists a permutation $\sigma$ of the coordinates of $x$, that sends $x$ to $y$. This happens
because the elements of $X$ are $n$-tuples obtained by permutation of the coordinates of the fixed $n$-tuple $(z_1, \dots, z_n)$. Suppose that
$\sigma=\sigma_{i_1} \dots \sigma_{i_l}$, then $\tau:=\tau_{i_1} \dots, \tau_{i_l}$ satisfies that $\phi(\tau)(v_x)=\lambda v_y$, for some
non-zero complex number $\lambda$. Then $W$ contains $v_y$ and therefore, $W$ contains the basis $\beta=\{v_x : x\in X\}$.

As $\phi(\tau_k)$ is a self-adjoint operator, it commutes with $P_W$, the orthogonal projection over the subspace $W$. Therefore,
$(\phi(\tau_k))^2$ commute with $P_W$. On the other hand, note that $(\phi(\tau_k))^2(v_x)=|q_{x_k, x_{k+1}}|^2 v_x$, hence, $(\phi(\tau_k))^2$
is diagonal in the basis $\beta=\{v_x : x\in X\}$. Then, the matrix of $P_W$ has at least the same blocks than $(\phi(\tau_k))^2$ for all $k$,
$1\leq k \leq n-1$.

If for some $k$, the matrix of $(\phi(\tau_k))^2$ has one block of size $1\times 1$, then the matrix of $P_W$ has one block of size $1 \times
1$. In other words, there exists $x\in X$ such that $v_x$ is an eigenvector. If the eigenvalue associated to $v_x$ is non-zero, then $v_x\in
W$.

It rest to see that the matrix of $(\phi(\tau_k))^2$ has all its blocks of size $1\times 1$. By hypothesis, for each pair of vectors in the
basis $\beta$,  $v_x$ and $v_y$, there exists $k$, $1\leq k \leq n-1$, such that $|q_{x_k, x_{k+1}}|^2\neq |q_{y_k, y_{k+1}}|^2$. Fix any order
in $X$ and let $x$ and $y$ the first and second element of $X$. Then there exists $k$ such that $v_x$ and $v_y$ are eigenvectors of
$(\phi(\tau_k))^2$ of different eigenvalue. Hence $(\phi(\tau_k))^2$ has the first block of size $1\times 1$. As $(\phi(\tau_j))^2$ commute with
$(\phi(\tau_k))^2$ for all $j$, $(\phi(\tau_j))^2$ also has this property.

By induction, suppose that for all $j$ $(\phi(\tau_j))^2$ has its $r-1$ first blocks of size $1\times 1$. Let $x', y'$ the elements $r$ and
$r+1$ of $X$, then there exists $k'$ such that $v_{x'}$ and $v_{y'}$ are eigenvectors of $(\phi(\tau_{k'}))^2$ of different eigenvalue. Hence,
$(\phi(\tau_{k'}))^2$ has the $r$ block of size $1\times 1$. Therefore $(\phi(\tau_j))^2$ too because it commute with $(\phi(\tau_{k'}))^2$, for
all $j$. Then we obtain that all the blocks are of size $1\times 1$.

\end{proof}

Note that if the numbers $q_{x_k, x_{k+1}}$ are all equal and $|X|>1$, then $\phi$ is not irreducible because the subspace $W$, generated by the
vector $v=\sum_{x\in X} v_x$, is an invariant subspace.

\subsection{Examples} We are going to compute some explicit examples of this family of representations. We will show that the standard
representation (\cite{S}, \cite{TYM}) is a member of this family.

\subsubsection{\it{Standard Representation}} Let $z_1=1$ and $z_j=0$ for all $j=2, \dots, n$. Then the cardinality of $X$ is $n$ and $\dim V=n$
too. For each $x\in X$, let $q_{x_k, x_{k+1}}=1 +(t-1) x_{k+1}$, where $t\neq 0,1$ is a complex number. Therefore $\phi:\B_n \rightarrow
\Aut(V)$, given by $\phi(\tau_k)v_x=q_{x_k, x_{k+1}} v_{\sigma_k(x)}$, is equivalent to the standard representation $\rho$, given by $$
\rho(\tau_k)= \left(\begin{array}{ccccccc}
          \  1      &\  \   &\  \      &\  \   &\  \   &\  \        &\  \  \\
          \  \      &\  1   &\  \      &\  \   &\  \   &\  \        &\  \  \\
          \  \      &\  \   &\  \ddots &\  \   &\  \   &\  \        &\  \  \\
          \  \      &\  \   &\  \      &\  0   &\  t   &\  \        &\  \  \\
          \  \      &\  \   &\  \      &\  1   &\  0   &\  \        &\  \  \\
          \  \      &\  \   &\  \      &\  \   &\  \   &\  \ddots   &\  \  \\
          \  \      &\  \   &\  \      &\  \   &\  \   &\  \        &\  1
\end{array}\right) $$ where $t$ is in the place $(k, k+1)$. In fact, if $\{\beta_j: j=1, \dots, n\}$ is the canonical basis of $\C^n$, and if
$x_j$ is the element of $X$ with $1$ in the place $j$ and zero elsewhere, define $$ \begin{array}{rcl}
 \alpha :\C^n & \rightarrow &V  \\
      \beta_j & \mapsto     &v_{x_j}
\end{array} $$

Then $\alpha (\rho(\tau_k) (\beta_j))=\phi(\tau_k)(\alpha(\beta_j))$ for all $j=1, \dots, n$. Hence the representations are equivalent.

\subsubsection{Example}

Let $z_1, \dots, z_n\in \{0, 1\}$, such that $z_1=z_2=\dots =z_m=1$ and $z_{m+1}=\dots =z_n=0$. Then the cardinality of $X$ is $\left(
\begin{smallmatrix} n\\\noalign{\medskip}m \end{smallmatrix}\right)=\frac{n!}{m! (n-m)!}$. If $V_m$ is the vector space with basis
$\beta_m=\{v_x : x\in X\}$, then $\dim V_m=\frac{n!}{m! (n-m)!}$.

For each $x:=(x_1, \dots, x_n)\in X$, let $$q_{x_k, x_{k+1}}=\left \{\begin{array}{ll}
    1       &\text{    if   } x_k=x_{k+1}   \\
    t       &\text{    if   } x_k\neq x_{k+1}
\end{array}\right. $$ where $t$ is a real number, $t\neq 0, 1, -1$.

We define $\phi_m: \B_n\rightarrow \Aut (V_m)$, given by $$ \phi_m(\tau_k)v_x=q_{x_k, x_{k+1}} v_{\sigma_k(x)} $$ For example, fixing the
lexicographic order in $X$, if $n=5$ and $m=3$, then $\dim V_m=10$, the ordered basis is $$ \begin{aligned} \beta:=\{&v_{(0,0,1,1,1)},
v_{(0,1,0,1,1)}, v_{(0,1,1,0,1)},v_{(0,1,1,1,0)}, v_{(1,0,0,1,1)}, \\
     &v_{(1,0,1,0,1)},v_{(1,0,1,1,0)},v_{(1,1,0,0,1)}, v_{(1,1,0,1,0)}, v_{(1,1,1,0,0)}\}
\end{aligned} $$ and the matrices in this basis are $$ \phi_3(\tau_1)= \left(\begin{array}{cccccccccc}
          \  1      &\  \   &\  \      &\  \   &\  \   &\  \        &\  \   &\  \   &\  \    &\  \  \\
          \  \      &\  0   &\  \      &\  \   &\  t   &\  \        &\  \   &\  \   &\  \    &\  \  \\
          \ \       &\  0   &\  0      &\  \   &\  0   &\  t        &\  \   &\  \   &\  \    &\  \  \\
          \ \       &\  0   &\  0      &\  0   &\  0   &\  0        &\  t   &\  \   &\  \    &\  \  \\
          \ \       &\  t   &\  0      &\  0   &\  0   &\  0        &\  0   &\  \   &\  \    &\  \  \\
          \ \       &\  \   &\  t      &\  0   &\  \   &\  0        &\  0   &\  \   &\  \    &\  \  \\
          \ \       &\  \   &\  \      &\  t   &\  \   &\  \        &\  0   &\  \   &\  \    &\  \  \\
          \  \      &\  \   &\  \      &\  \   &\  \   &\  \        &\  \   &\  1   &\  \    &\  \  \\
          \  \      &\  \   &\  \      &\  \   &\  \   &\  \        &\  \   &\  \   &\  1    &\  \  \\
          \  \      &\  \   &\  \      &\  \   &\  \   &\  \        &\  \   &\  \   &\  \    &\  1
\end{array}\right) $$ $$ \phi_3(\tau_2)= \left(\begin{array}{cccccccccc}
          \  0      &\  t   &\  \      &\  \   &\  \   &\  \        &\  \   &\  \   &\  \    &\  \  \\
          \  t      &\  0   &\  \      &\  \   &\  \   &\  \        &\  \   &\  \   &\  \    &\  \  \\
          \  \      &\  \   &\  1      &\  \   &\  \   &\  \        &\  \   &\  \   &\  \    &\  \  \\
          \  \      &\  \   &\  \      &\  1   &\  \   &\  \        &\  \   &\  \   &\  \    &\  \  \\
          \  \      &\  \   &\  \      &\  \   &\  1   &\  \        &\  \   &\  \   &\  \    &\  \  \\
          \  \      &\  \   &\  \      &\  \   &\  \   &\  0        &\  \   &\  t   &\  \    &\  \  \\
          \  \      &\  \   &\  \      &\  \   &\  \   &\  0        &\  0   &\  0   &\  t    &\  \  \\
          \  \      &\  \   &\  \      &\  \   &\  \   &\  t        &\  0   &\  0   &\  0    &\  \  \\
          \  \      &\  \   &\  \      &\  \   &\  \   &\  \        &\  t   &\  \   &\  0    &\  \  \\
          \  \      &\  \   &\  \      &\  \   &\  \   &\  \        &\  \   &\  \   &\  \    &\  1
\end{array}\right) $$ $$ \phi_3(\tau_3)= \left(\begin{array}{cccccccccc}
          \ 1      &\ \   &\ \      &\ \   &\ \   &\ \        &\ \   &\ \   &\ \     &\ \   \\
          \ \      &\ 0   &\ t      &\ \   &\ \   &\  \       &\ \   &\ \   &\ \     &\ \   \\
          \  \     &\ t   &\ 0      &\ \   &\ \   &\ \        &\ \   &\ \   &\  \    &\ \   \\
          \  \     &\ \   &\ \      &1 \   &\ \   &\ \        &\ \   &\ \   &\ \     &\ \   \\
          \  \     &\ \   &\ \      &\ \   &\ 0   &\ t        &\ \   &\ \   &\ \     &\ \   \\
          \  \     &\ \   &\ \      &\ \   &\ t   &\ 0        &\ \   &\ \   &\ \     &\ \   \\
          \  \     &\  \  &\ \      &\ \   &\ \   &\ \        &\ 1   &\ \   &\ \     &\ \    \\
          \  \     &\  \  &\ \      &\ \   &\ \   &\ \        &\ \   &\ 1   &\ \     &\ \    \\
          \  \     &\  \  &\ \      &\ \   &\ \   &\ \        &\ \   &\ \   &\  0    &\  t  \\
          \  \     &\ \   &\ \      &\ \   &\ \   &\ \        &\ \   &\ \   &\  t    &\  0
\end{array}\right) $$ $$ \phi_3(\tau_4)= \left(\begin{array}{cccccccccc}
          \  1      &\  \   &\  \      &\  \   &\  \   &\  \        &\  \   &\   &\  \    &\  \  \\
          \  \      &\  1   &\  \      &\  \   &\  \   &\  \        &\  \   &\   &\  \    &\  \  \\
          \  \      &\  \   &\  0      &\  t   &\  \   &\  \        &\  \   &\   &\  \    &\  \  \\
          \  \      &\  \   &\  t      &\  0   &\  \   &\  \        &\  \   &\   &\  \    &\  \  \\
          \  \      &\  \   &\  \      &\  \   &\  1   &\  \        &\  \   &\   &\  \    &\  \  \\
          \  \      &\  \   &\  \      &\  \   &\  \   &\  0        &\ \  t   &\   &\  \    &\  \  \\
          \  \      &\  \   &\  \      &\  \   &\  \   &\  t        &\ \  0   &\   &\  \    &\  \  \\
          \  \      &\  \   &\  \      &\  \   &\  \   &\  \        &\  \   &0   &\  t    &\  \  \\
          \  \      &\  \   &\  \      &\  \   &\  \   &\  \        &\  \   &t   &\  0    &\  \  \\
          \  \      &\  \   &\  \      &\  \   &\  \   &\  \        &\  \   &\  \   &\    &\  1
\end{array}\right) $$

With this notation, we have the following results,

\begin{theorem} Let $n>2$, then $(\phi_m, V_m)$ is an irreducible representations of $\B_n$, for all $1\leq m<n$ such that $2m\neq n$.

If $n=2m$ then $(\phi_m, V_m)$ is sum of two irreducible representations of $\B_n$.
 \end{theorem} \begin{proof} Suppose that $n\neq 2m$. Let $x\neq y\in X$, then there exists $j$, $1\leq j\leq n$, such that
$x_j\neq y_j$. If $j>1$, we may suppose that $x_{j-1}=y_{j-1}$,  then $q_{x_{j-1}, x_j}\neq q_{y_{j-1}, y_j}$, therefore $|q_{x_{j-1},
x_j}|^2\neq |q_{y_{j-1}, y_j}|^2$. If $j=1$, and $n\neq 2m$, there exists $l=2, \dots, n$ such that $x_{l-1}\neq y_{l-1}$ and $x_l=y_l$,  then
$|q_{x_{l-1}, x_l}|^2\neq |q_{y_{l-1}, y_l}|^2$. Then, by theorem \ref{irreducible}, $\phi_m$ is an irreducible representation.

Note that if $n=2m$, then $X=\{x, y_x : x\in Y\}$, for some $Y\subset X$, where $y_x$ is obtained of $x$ changed the zeros by ones and the ones by zeros. For example if
$$
X:=\{(0,0,1,1), (0,1,0,1), (0,1,1,0), (1,0,0,1),(1,0,1,0), (1,1,0,0) \}
$$
then $Y:=\{(0,0,1,1), (0,1,0,1), (0,1, 1, 0)\}$.

Let $W_1, W_2$ be the vector spaces generated by $\beta_1:=\{v_x+v_{y_x}: x\in Y\}$ and $\beta_2:=\{v_x-v_{y_x}: x\in Y\}$ respectively, then $W_1$ and $W_2$ are irreducible invariant subspaces of $V_m$.
 \end{proof}

The {\it{corank }} of a finite dimensional representation $\phi$ of $\B_n$ is the rank of $(\phi(\tau_k)-1)$. This number does not depend on $k$
because all the $\tau_k$ are conjugate to each other (see p. 655 of \cite{C}).

\begin{theorem} If $n>2$, $1\leq m <n$ and $2m\neq n$ then $(\phi_m, V_m)$ is an irreducible representation of dimension $\left( \begin{smallmatrix}
n\\\noalign{\medskip}m \end{smallmatrix}\right)$ and corank $\frac{2(n-2)!}{(m-1)!(n-m-1)!}$. \end{theorem} \begin{proof} By theorem before,
$(\phi_m, V_m)$ is an irreducible representation. The dimension of $\phi$ is the cardinality of $X$, then

$$ \dim V_m= \left( \begin{smallmatrix} n\\\noalign{\medskip}m \end{smallmatrix}\right)=\frac{n!}{m! (n-m)!} $$

 We compute the corank of $\phi_m$. Let $x\in X$ such that $\sigma_k(x)= x$, then $x_k=x_{k+1}$ and $q_{x_k, x_{k+1}}=1$.
Therefore $\phi_m(\tau_k)(v_x) =v_x$. Hence the corank of $\phi_m$ is equal to the cardinality of $Y=\{x\in X : \sigma_k(x)\neq x \}$. But it is
equal to the cardinality of $X$ minus the cardinality of $\{x\in X : x_k=x_{k+1}=0 \text{  or   } x_k=x_{k+1}=1 \}$. Therefore $$
\begin{aligned} cork(\phi_m)=rk(\phi_m(\tau_k) -1)&=\frac{n!}{m!(n-m)!}-\frac{(n-2)!}{m! (n-m-2)!}-\frac{(n-2)!}{(m-2)!(n-m)!}\\
                                  &=\frac{2(n-2)!}{(m-1)!(n-m-1)!}
\end{aligned} $$ \end{proof}

In the example $n=5$ and $m=3$, we have that $cork(\phi_m)=6$.

Note that if $m=1$, the dimension of $\phi_m$ is $n$ and the corank is $2$. Therefore $\phi_1$ is equivalent to the
 standard representation, because this is the unique irreducible representations of $\B_n$ of dimension $n$ \cite{S}.

\section*{Acknowledgment} The authors thanks to Aroldo Kaplan for his helpful comments. Also, the authors thanks to Vaughan Jones for pointed out the mistake in the previous version.

\end{document}